\documentclass[a4paper,11pt]{amsart}
\usepackage{graphicx}
\usepackage{psfrag}
\usepackage[centering]{geometry}
\usepackage{amssymb}
\usepackage[colorlinks=true]{hyperref}
\gdef\s{s}
\title[On the doubly refined enumeration of ASM{\s} and TSSCPP{\s}]{On the Doubly Refined Enumeration of Alternating Sign Matrices and Totally Symmetric Self-Complementary Plane Partitions}
\author[T.~Fonseca]{Tiago~Fonseca}
\address{Tiago Fonseca, LPTHE (CNRS, UMR 7589), Univ Pierre et Marie Curie-Paris6, 75252 Paris Cedex, France.}
\email{fonseca\,@\,lpthe.jussieu.fr}
\author[P.~Zinn-Justin]{Paul~Zinn-Justin}
\address{Paul Zinn-Justin, LPTMS (CNRS, UMR 8626), Univ Paris-Sud, 91405 Orsay Cedex, France; and LPTHE (CNRS, UMR 7589), Univ Pierre et Marie Curie-Paris6, 75252 Paris Cedex, France.}
\email{pzinn\,@\,lpthe.jussieu.fr}
\thanks{PZJ was supported by
EU Marie Curie RTN %Research Training Networks
``ENRAGE'' MRTN-CT-2004-005616, ``ENIGMA'' MRT-CT-2004-5652,
ESF program ``MISGAM''
and ANR program ``GIMP'' ANR-05-BLAN-0029-01.}
\thanks{The authors thank N.~Kitanine for discussions, and J.-B.~Zuber
for a careful reading of the manuscript.}
\numberwithin{equation}{section}
\newtheorem*{conj*}{Theorem}
\newtheorem{lemma}{Lemma}

\date{March 2008}

\begin{document}

\begin{abstract}
We prove the equality of 
doubly refined enumerations of Alternating Sign Matrices
and of Totally Symmetric Self-Complementary Plane Partitions
using integral formulae originating from certain solutions
of quantum Knizhnik--Zamolodchikov \hbox{equation}.
\end{abstract}

\maketitle
{\footnotesize\tableofcontents}

\section{Introduction}
It is the purpose of this work to revisit an old problem using some new ideas.
The old problem is the interconnection between two distinct classes of combinatorial objects whose enumerative properties are intimately related: Alternating
Sign Matrices and Plane Partitions \cite{Bressoud}. The new ideas come from recent developments in the so-called Razumov--Stroganov conjecture (formulated in
\cite{RSc}; see also \cite{XXZ-ASM-PP,dG-02}). The
Razumov--Stroganov conjecture identifies the entries of the 
Perron--Frobenius vector of a certain stochastic matrix with cardinalities
of subsets of Alternating Sign Matrices, the latter being
reinterpreted as configurations of a certain two-dimensional
statistical model (so-called Fully Packed Loops). Even though this statement
is still a conjecture, some progress has been made in this area in a
series of papers by Di Francesco and Zinn-Justin, starting with
\cite{DFZJ-sr}. The method they used was, as it turned out,
equivalent to finding appropriate polynomial solutions of the quantum 
Knizhnik--Zamolodchikov equation \cite{DFZJ-qKZ}. Integral representations 
for these and their relation to plane partition enumeration 
were discussed in \cite{DFZJ-qKZ-TSSCPP}; we shall use these
integral formulae in the present work (noting that 
these can be considered as purely formal integrals,
so they are simply a way of encoding generating functions).

The paper is organized as follows. In section 2, we define the various
combinatorial objects and corresponding statistical models that will
be needed. In section 3, we formulate the main theorem of the paper:
the equality of doubly refined enumerations of Alternating Sign Matrices
and of Totally Symmetric Self-Complementary Plane Partitions. Section 4
contains the proof, based on the use of integral formulae. Finally,
the appendices contain various technical results that are needed in the proof.
Note that even though we use some concepts and methods from
exactly solvable statistical models, this paper is self-contained and
all proofs are purely combinatorial in nature.

\section{The models}

In this section we define the various models that appear in this work. There are two distinct models. On the one hand we have Alternating Sign Matrices (ASMs) which are in bijection with configurations of the 6-Vertex model (also known as ice model) with Domain Wall Boundary Conditions, as well as with Fully Packed Loop configurations (FPL). Here we only discuss ASMs and 6-V model.

On the other hand we have Totally Symmetric Self-Complementary Plane Partitions, which are in bijection with a certain class of Non-Intersecting Lattice Paths.

\subsection{Alternating Sign Matrices} \label{ASM}

An Alternating Sign Matrix (ASM) is a square matrix made of 0s, 1s and -1s 
such that if one ignores 0s,
1s and -1s alternate on each row and column starting and ending with 1s.
Here are all $3\times 3$ ASMs:

\[
\begin{smallmatrix}1&0&0\\0&1&0\\0&0&1\end{smallmatrix}
\qquad
\begin{smallmatrix}1&0&0\\0&0&1\\0&1&0\end{smallmatrix}
\qquad
\begin{smallmatrix}0&1&0\\1&0&0\\0&0&1\end{smallmatrix}
\qquad
\begin{smallmatrix}0&1&0\\0&0&1\\1&0&0\end{smallmatrix}
\qquad
\begin{smallmatrix}0&1&0\\1&\!-1\!&1\\0&1&0\end{smallmatrix}
\qquad
\begin{smallmatrix}0&0&1\\1&0&0\\0&1&0\end{smallmatrix}
\qquad
\begin{smallmatrix}0&0&1\\0&1&0\\1&0&0\end{smallmatrix}
\]

Thus, there are exactly $7$ ASMs of size $n=3$. 

These matrices have been studied by Mills, Robbins and Rumsey since the early 1980s \cite{MRR-Mac, MRR-ASM, RR-ASM, MRR-TSSCPP}. It was then
conjectured that $A_n$, the number of ASMs of size $n$, is given by:
\begin{equation}
A_n = \prod_{j=0}^{n-1} \frac{(3 j +1)!}{(n+j)!}=1,\ 2,\ 7,\ 42,\ 429,\ \ldots
\end{equation}

This was subsequently proved by Zeilberger in 1996 in an 84 page article \cite{Zeil-ASM}. A shorter proof was given by Kuperberg \cite{kup-ASM} in 1998.
The latter is based on the equivalence to the 6-V model, which we shall also
use here.

\subsection{6-Vertex model}

Let us now turn to the 6-Vertex Model. 
The model consists in a square grid of size $n\times n$ in which each edge is given an orientation (an arrow), such that at each vertex there are two arrows pointing in and two arrows pointing out. We use here some very specific boundary conditions (Domain Wall Boundary Conditions, DWBC):
all arrows at the left and the right are pointing in and at the bottom and the top are pointing out.

\begin{figure}

\centering
\includegraphics[scale=0.4]{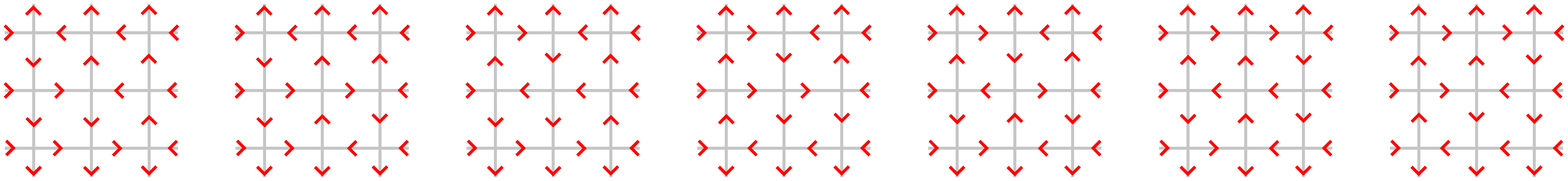}
\caption{\label{6V3}The 6-Vertex Model is defined on a $n\times n$ grid. To each link in the network we associate an arrow which can take two directions, the only constraint being that at each site there are two arrows pointing in and two arrows pointing out (this leaves 6 possible vertex configurations). We are only interested in the configurations such that the arrows at the top and at the bottom are pointing out and the arrows at the left and the right are pointing in. Here we draw all states possibles for $n=3$.}

\end{figure}

\begin{figure}

\psfrag{a}[0][0][1][0]{0}\psfrag{b}[0][0][1][0]{1}\psfrag{c}[0][0][1][0]{-1}
\centering
\includegraphics[scale=0.8]{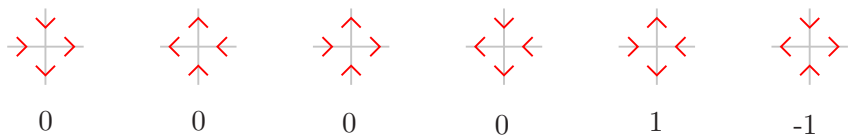}
\caption{\label{6VASM}Rules to replace each vertex of a 6-V configuration with a $0$ or $\pm1$. Conversely, one can consistently build a 6-V configuration from an ASM starting from the fixed arrows on the boundary, continuing arrows through the $0s$ and reversing them through the $\pm 1$.}
\end{figure}

On figure~\ref{6V3} we draw all the possible configurations at $n=3$.  There are once again 7 configurations of size $n=3$. Indeed, there is an easy bijection between ASMs and 6-V configurations with DWBC, which is described schematically on figure~\ref{6VASM}.

\subsection{Totally Symmetric Self-Complementary Plane Partitions}

We describe here Plane Partitions in two different ways, either pictorially or as arrays of numbers.

\begin{figure}
 
\centering
\includegraphics[scale=0.5]{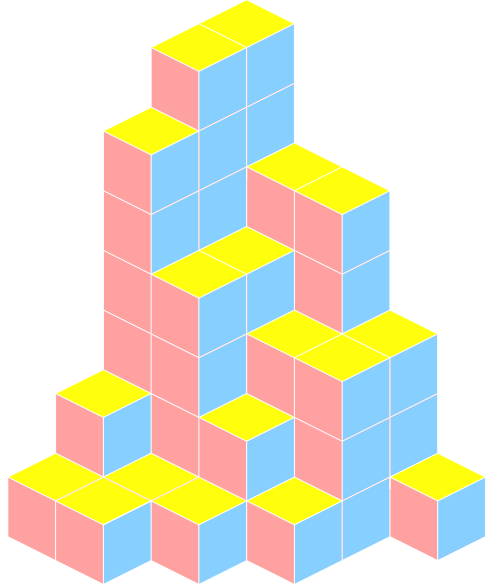}
\caption{\label{PP}We can see a plane partition (PP) as a stack of unit cubes pushed into a corner.}
\end{figure}

Pictorially, a plane partition is a stack of unit cubes pushed into a corner (gravity pushing them to the corner) and drawn in isometric perspective, as examplified on figure~\ref{PP}.

An equivalent way of describing these objects is to form the array of heights of each stack of cubes. In this formulation the effect of ``gravity'' is that each number in the array is less or equal than the numbers immediately above and to the left. For example the plane partition on figure~\ref{PP} may be translated into the array

\[\begin{array}{l} 75531\\ 7433 \\6421 \\211 \\11 \end{array}\]

Plane partitions were first introduced by MacMahon in 1897.
A problem of interest is the enumeration of plane partitions that have some specific symmetries. The Totally Symmetric and Self-Complementary Plane Partitions (TSSCPPs) are one of these symmetry classes. In the pictorial representation,
they are Plane Partitions inside a $2n\times 2n\times 2n$ cube which are invariant under the following symmetries:
all permutations of the axes of the cube of size $2n\times 2n\times 2n$; and taking the complement, that is putting cubes where they are absent and vice versa, and flipping the resulting set of cubes to form again a Plane Partition.

Alternatively, they can be described as $2n\times 2n$ arrays of heights.
In the $n=3$ case, we have, once again, $7$ possible configurations:
\begin{equation}
\begin{array}{ccccccc}
666333 & 666433 & 666433 & 666543 & 666543 & 666553 & 666553\\
666333 & 666333 & 666433 & 665332 & 665432 & 655331 & 655431\\
666333 & 665332 & 664322 & 655331 & 654321 & 655331 & 654321\\
333000 & 433100 & 443200 & 533110 & 543210 & 533110 & 543210\\
333000 & 333000 & 332000 & 433100 & 432100 & 533110 & 532110\\
333000 & 332000 & 332000 & 321000 & 321000 & 311000 & 311000\end{array} \label{TSSCPP}
\end{equation}
and more generally we obtain $A_n$ for any $n$. In fact Zeilberger's proof
of the ASM conjecture amounts to showing (non-bijectively) that
ASMs and TSSCPPs are equinumerous.

\subsection{Non-Intersecting Lattice Paths} \label{NILP sec}

Another important class of objects consists of Non-Intersecting Lattice Paths (NILPs). 
These paths are defined in a lattice and connect a set of initial points to a set of final points following certain rules (see Ref.~\cite{Lind,GV} for the general framework). The most important feature of NILPs is that the various paths do not touch one another.

In order to better understand the bijection between NILPs and TSSCPPs,
it is convenient to consider an intermediate class of objects: Non-Crossing Lattice Paths (NCLPs), which are similar to NILPs except for the fact the paths are allowed to share a common site, although they are still forbidden to cross each other.

We proceed with the description of the bijection between TSSCPPs and a class of NCLPs. Each TSSCPP is defined by a subset of numbers of the arrays of~\eqref{TSSCPP}, a possible choice is the triangles at the bottom right:
\[\begin{array}{lllllll}
0 & 1 & 2 & 1 & 2 & 1 & 2\\
00 & 00 & 00 & 10 & 10 & 11 & 11\\
000 & 000 & 000 & 000 & 000 & 000 & 000\end{array}\]
It is easy to prove that this part of the array together with the symmetries which characterize the TSSCPPs are enough to reconstruct the whole TSSCPP.

\begin{figure}

\centering
\includegraphics[scale=0.8]{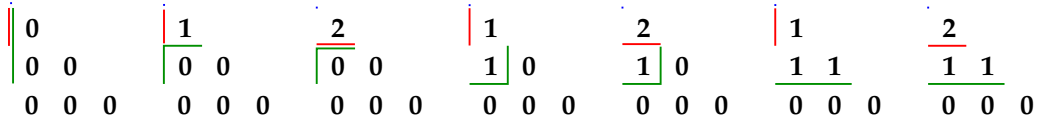}
\caption{Reformulation of TSSCPPs as NCLPs, in the example of size $n=3$. 
If the origin is at the upper right corner, then at each point $(0,-i)$, $i\in \{0,\ldots,n-1\}$, begins a path which can only go upwards or to the right, and stops when it reaches the diagonal $(j,-j)$, in such a way that
the numbers below/to the right of it are exactly those 
less or equal to $n-i$.\label{NILPX}}

\end{figure}

Then, we draw paths separating the different numbers appearing, as explained on figure~\ref{NILPX}.

\begin{figure}

\centering
\includegraphics[scale=0.6]{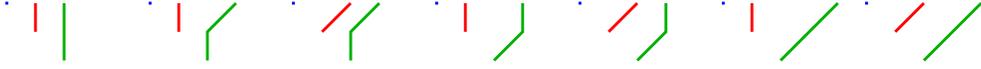}
\caption{We transform our NCLPs into NILPs: the starting point is now shifted to the right, and the horizontals steps become diagonal steps. %Avoiding that the paths touch one another.
\label{NILP} }

\end{figure}

The bijection with the NILPs is easily achieved by shifting the paths (NCLPs) according to the following rules:

\begin{itemize}

 \item The $i^{\textrm{th}}$ path begins at $(i,-i)$;

 \item The vertical steps are conserved and the horizontal steps ($\rightarrow$) are replaced by diagonal steps ($\nearrow$).

\end{itemize}
An example ($n=3$) is shown on figure~\ref{NILP}.

\begin{figure}

\centering
\includegraphics[scale=0.63]{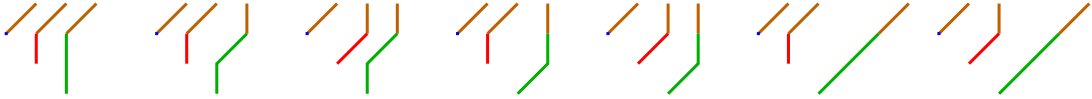}
\caption{To each path we add one extra step in order that two final points consecutive differ by an odd number. The first extra step is diagonal.\label{NILPplus}}

\end{figure}

Our last modification is the addition of one extra step to all paths. To the first path we add a diagonal step, as for the other paths the choice is made such that the difference between the final point of two consecutive paths is an odd number, as examplified on figure~\ref{NILPplus}.

\section{The conjecture}

Various conjectures have been made to connect ASMs and TSSCPPs. Building on the already mentioned ASM conjecture by Mills and Robbins, which says that the number of ASMs of size $n$ is equal to the number of TSSCPPs of size $2n$ (and which is now a theorem), there are conjectures about ``refined'' enumeration. Before describing them we need some more definitions. 

\subsection{ASM generating function} \label{A-n}

Each ASM, as can be easily proven, has one and only one $1$ on the first row and on the last row. It is natural to classify ASMs according to its position. Therefore, we count the ASMs of size $n$ with the first $1$ in the $i^{\rm th}$ position and the last $1$ in the $j^{\rm th}$ position: $\tilde{A}_{n,i,j}$.

We build the corresponding generating function:
\begin{equation}
 \tilde{A}_n(x,y):=\sum_{i,j}\tilde{A}_{n,i,j}x^{i-1}y^{j-1}
\end{equation}

We define also $A_{n,i,j}$, which counts the ASMs with the first $1$ in the 
$i^{\rm th}$ column and the last $1$ in the $(n-j+1)^{\rm st}$ column:
\begin{equation}
 A_{n,i,j}=\tilde{A}_{n,i,n-j+1}
\end{equation}
And its generating function:
\begin{equation}
 A_n(x,y):=\sum_{i,j}A_{n,i,j}x^{i-1}y^{j-1}
\end{equation}

\subsubsection*{Some trivial symmetries}
By reflecting the ASMs horizontally and vertically one gets:
\[A_{n,i,j}=A_{n,j,i}\]
whereas by reflecting them only horizontally one gets:
\[A_{n,i,j}=A_{n,n-i+1,n-j+1}\]
Obviously these symmetries are also valid for $\tilde{A}_{n,i,j}$.

\subsection{NILP generating function} \label{N-X}

First we recall the definition of the type of NILPs used in this article, of size $n$:

\begin{itemize}
 \item The paths are defined on the square grid. Each step connects a site to a neighbor and can be either vertical (up $\uparrow$) or diagonal (up right $\nearrow$).

 \item There are $n$ starting points with coordinates $(i,-i),\ i\in \{0,1,..,n-1\}$. The endpoints are at $(i,0)$ (so that the length of the $i^{\rm th}$ path is $i$). 

 \item Paths do not touch each other.
\end{itemize}

It is convenient to add an extra step, as explained in section~\ref{NILP sec}, 
defined uniquely by the following:

\begin{itemize}
 \item Two consecutive paths, after the extra step, differ by an odd number.

 \item The extra step for the first path (at $(0,0)$) is diagonal.
\end{itemize}

Let $\alpha$ be a NILP, we define $u^0_n(\alpha)$ as the number of vertical steps in the extra step and $u^1_n(\alpha)$ as the number of vertical step in the last step of each path (see appendix~\ref{U-u} for an extended definition).

The generating function is:
\begin{equation}
 U^{0,1}_n(x,y):= \sum_\alpha x^{u^0(\alpha)} y^{u^1(\alpha)}= \sum_{i,j} U^{0,1}_{n,i,j} x^i y^j
\end{equation}
where $U^{0,1}_{n,i,j}$ is the number of NILPs of size $n$ with $i$ vertical extra steps and $j$ vertical last steps.

\subsection{The conjecture} \label{conj}
We now present the conjecture, formulated by Mills, Robbins and Rumsey in
a slightly different language (see section~\ref{Mills2} for a detailed
translation), whose proof is the main focus of the present work:

\begin{conj*}
 The number of ASMs of size $n$ with the $1$ of the first row in the $(i+1)^{\rm st}$ position and the $1$ of last row in the $(j+1)^{\rm st}$ position is the same as the number of NILPs (corresponding to TSSCPPs, and with the extra step) with $i$ vertical extra steps and $j$ vertical steps in the last step. Equivalently,
\[\tilde{A}_n(x,y)=U^{0,1}_n(x,y)\]
\end{conj*}

For example, at $n=3$, using the ASMs given in section~\ref{ASM} and the TSSCPPs given on figure~\ref{NILPplus}, we compute:
\begin{align*}
 \tilde{A}_3(x,y) = & y^2+y+xy^2+x+xy+x^2y+x^2\\
 U^{0,1}_n(x,y) = & y^2+xy+x^2+xy^2+x^2y+y+x
\end{align*}
This is the doubly refined enumeration.
Of course, by specializing one variable, one recovers
the simple refined enumeration, i.e.\ that
the number of ASMs of size $n$ with the $1$ of the first row in the $i+1$ position is the same as the number of NILPs (corresponding to the TSSCPPs and with the extra step) with $i$ vertical extra steps:
\[A_n(x):=\tilde{A}_n(x,1)=U^{0,1}_n(x,1):= U^{0}_n(x)\]
and by specializing two variables, that
the number of ASMs of size $n$ is the same as the number of TSSCPPs of size $2n$:
\[A_n=A_n(1)=U^{0}_n(1)\]

\section{The proof}

\subsection{ASM counting as the partition function of the 6-Vertex model} \label{6-V pf}

\begin{figure}
\centering
\psfrag{a}[0][0][1][0]{$a=q^{-1/2}w-q^{1/2}z$}
\psfrag{b}[0][0][1][0]{$b=q^{-1/2}z-q^{1/2}w$}
\psfrag{c}[0][0][1][0]{$c=(q^{-1}-q)z^{1/2}w^{1/2}$}
\includegraphics[scale=0.8]{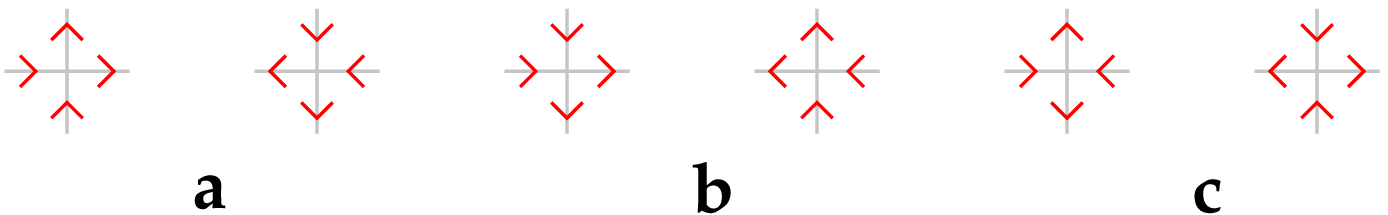}
\caption{To each site configuration corresponds a statistical weight. These weights depend on three parameters: $w$ (resp.\ $z$) which characterizes the column (resp.\ row), and a global parameter $q$ which will be eventually specialized to a cubic root of unity.\label{poids}}

\end{figure}

In order to solve the ASM enumeration problem, it is convenient to generalize it by considering weighted enumeration. This amounts to computing the partition function of the 6-Vertex model, that is the summation over 6-V configurations with DWBC such that to each vertex is given a statistical weight, as shown on figure~\ref{poids}, depending on $n$ horizontal spectral parameters (one for each row) $\{z_{1},z_{2},\ldots,z_{n}\}$, $n$ vertical spectral parameters $\{z_{n+1},z_{n+2},\ldots,z_{2n}\}$ and one global parameter $q$. This computation was performed by Izergin \cite{Iz-6V}, using recursion relations written by Korepin \cite{kor}, and the result is a $n\times n$ determinant (IK determinant). It is a symmetric function of the set $\{z_1,\ldots z_n\}$ and of the set of $\{z_{n+1},\ldots,z_{2n}\}$. Much later, it was observed by Stroganov \cite{strog}
and Okada \cite{okada} 
that when $q = e^{2 \pi i /3}$, the partition function is totally symmetric, 
i.e.\ in the full set $\{z_1,\ldots,z_{2n}\}$.
 
More precisely, if we denote by $\tilde Z_n$ the partition function,
and
\[
Z_n=
(-1)^{n(n-1)/2}{(q^{-1} - q)}^{-n} \prod _{i=1} ^{2n} z_i^{-1/2}
\ \tilde Z_n
\]
then $Z_n$ was identified with the Schur function corresponding to the Young diagram $Y_n$ with two rows of length $n-1$, two rows of length $n-2$, \dots, two rows of length 2 and two rows of length 1:
\begin{equation} 
Z_n(z_1,\ldots\ z_{2n})=s_{Y_n}(z_1,\ldots\ ,z_{2n})=\frac{\det[z_i ^{2n-j+d_j}]}{\det[z_i ^{2n-j}]} \label{schur}
\end{equation}
where $d_j$ is the sequence $\{n-1,n-1,n-2,n-2,\ldots,2,2,1,1,0,0\}$.
This formula is proved in appendix ~\ref{Z=Z' Z}, 
though its explicit form will not
be needed in what follows.

With this method we recover the unweighted enumeration by setting all $z_i=1$:
\begin{equation}
3^{-n(n-1)/2} Z_n(1,\ldots,1)=A_n
\end{equation}
where we recall that
$A_n$ is the number of ASMs of size $n$ (as explained in \ref{ASM}).

The case of interest to us is when all $z_i=1$ except $z_1$ and $z_{2n}$:
\begin{align*}
z_1  = & \frac{1+qt}{q+t}\\
z_{2n}  = & \frac{1+qu}{q+u}
\end{align*}
Using the fact
that $Z_n(z_1,\ldots,z_{2n})$ is a symmetric function of its arguments
(see appendix~\ref{Z=Z' Z}), we have
\[
Z_n(z_1=\frac{1+qt}{q+t},1,\ldots,1,z_{2n}=\frac{1+qu}{q+u})=Z_n(z_1=\frac{1+qt}{q+t},1,\ldots,1,z_n=\frac{1+qu}{q+u},1,\ldots,1)
\] 

The corresponding weights take the form
\begin{align*}
a_{x} & =  q^{-\frac{1}{2}}-q^{\frac{1}{2}}\left(\frac{1+qx}{q+x}\right)=\frac{q^{\frac{1}{2}}x}{q+x}(q^{-1}-q)\\
b_{x} & =  q^{-\frac{1}{2}}\left(\frac{1+qx}{q+x}\right)-q^{\frac{1}{2}}=\frac{q^{\frac{1}{2}}}{q+x}(q^{-1}-q)\\
c_{x} & =  (q^{-1}-q)\sqrt{\frac{1+qx}{q+x}}
\end{align*}

The partition function $\tilde Z_n$ becomes
\begin{align*}
\tilde{Z_{n}}  &=  (-i\sqrt{3})^{n^{2}-2n}\sum_{j,k}a_{t}^{j-1}b_{t}^{n-j}c_{t}a_{u}^{k-1}b_{u}^{n-k}c_{u}A_{n,j,k}\\
 & =  (-i\sqrt{3})^{n^{2}}\sqrt{\frac{1+qt}{q+t}}\sqrt{\frac{1+qu}{q+u}}\left(\frac{1}{q+t}\right)^{n-1}\left(\frac{1}{q+u}\right)^{n-1}q^{n-1}\sum_{j,k}t^{j-1}u^{k-1}A_{n,j,k}
\end{align*}
where $A_{n,j,k}$ is the number of ASMs of size $n$ such that the only $1$ in the first row is in column $j$ and the only $1$ in the last row is in column $n-k+1$. 

The normalization factor is equal to:
$(-1)^{n(n-1)/2}(-i\sqrt{3})^{n}\sqrt{\frac{1+qt}{q+t}}\sqrt{\frac{1+qu}{q+u}}$,
so we can finally compute
\begin{equation}
\frac{(q^{2}(q+t)(q+u))^{n-1}}{3^{n(n-1)/2}}Z_{n}=\sum_{j,k}t^{j-1}u^{k-1}A_{n,j,k}=A_n(t,u)
\end{equation}

Note that if one uses instead $z_{2n} = \frac{q+u}{1+qu}$, one gets the same formula, but with one index reversed
\begin{equation}
\frac{(q^{2}(q+t)(1+qu))^{n-1}}{3^{n(n-1)/2}}Z_{n}=\tilde{A}_n(t,u)
\end{equation}

\subsection{Integral formula for refined ASM counting}
The traditional expression for the partition function of the 6-V model
is the already mentioned IK formula. We shall not use it here.
We shall only need the following facts (true at $q=e^{2\pi i/3}$):

\begin{itemize}
\item $Z_1=1$.
\item $Z_n(z_1,\ldots,z_{2n})$ is a polynomial of degree $n-1$ in each variable.
%which is symmetric in the full set of variables
%$\{ z_1,\ldots,z_n\}$ and $\{ z_{n+1},\ldots,z_{2n}\}$.
%$\{ z_1,\ldots,z_{2n}\}$.
\item The $Z_n$ satisfy the recursion relation for all $i,j=1,\ldots,2n$
\begin{equation}\label{strogrecur}
Z_{n}(z_{1},\ldots,z_{j}=q^{2}z_{i},\ldots,z_{2n})=\prod_{k\ne i,j}(q z_{i}-z_{k})Z_{n-1}(z_1,\ldots,\hat z_i,\ldots,\hat z_j,\ldots,z_{2n})
\end{equation}
\end{itemize}
We recall how to prove them 
in appendix~\ref{Z=Z' Z} for the sake of completeness.

Furthermore, we need the following lemma
\begin{lemma}\label{wheel} A polynomial $P$ of degree $n-1$ in each variable
$z_1,\ldots,z_{2n}$ which satisfies the ``wheel condition''
\[
P(\ldots,z_i=z,\ldots,z_j=q^2 z,\ldots,z_k=q^4 z,\ldots)=0 \qquad \textrm{for\ all\ } i<j<k
\]
is entirely determined by its $c_n:=(2n)!/n!/(n+1)!$ values
at the following specializations: $(q^{\epsilon_1},\ldots,q^{\epsilon_{2n}})$ for all possible choices of $\{\epsilon_i=\pm 1\}$ such that $\sum_{i=1}^{2n} \epsilon_i=0 $ and $\sum_{i=1}^j\epsilon_i\leq 0$ for all $j\leq2n$.
\end{lemma}
This lemma is proved in appendix~\ref{Z=Z' lemma}.

The strategy is now to introduce a certain integral representation of the
partition function of the 6-V model with DWBC, say $Z'_n$
\begin{multline}
Z'_n:=(-1)^{\binom{n}{2}}\prod_{i<j}^{2n}(qz_{i}-q^{-1}z_{j})\label{Z'n}\\
\oint\ldots\oint\prod_{l}^{n}\frac{dw_{l}}{2\pi i}\frac{(qz_{2l-1}-q^{-1}w_l)\prod_{l<m}(w_{m}-w_{l})(qw_{l}-q^{-1}w_{m})}{\prod_{i\leq2l-1}(w_{l}-z_{i})\prod_{i\ge2l-1}(qw_{l}-q^{-1}z_{i})}
\end{multline}
where the integration contours surround counterclockwise
the $z_i$ (but not the $q^{-2}z_i$), 
and to show that $Z_n$ and $Z'_n$ are both polynomials of 
degree $n-1$ in each variable which  satisfy the ``wheel condition''
and coincide at the $c_n$ specializations of lemma~\ref{wheel}.

Let us first check that $Z_n$ satisfies the wheel condition. This is a direct
consequence of Eq.~(\ref{strogrecur}) in which one sets $z_k=q^4 z_i$.
It is equally 
straightforward to calculate $Z_n$ at the $c_n$ points of the lemma.
The computation goes inductively using Eq.~(\ref{strogrecur}) and it is left to the reader to check that
\[
Z_n(q^{\epsilon_1},\ldots,q^{\epsilon_{2n}})=3^{\binom{n}{2}}
\]

We now show that $Z'_n$ also satisfies the hypotheses of the lemma.
We proceed in steps.

\subsubsection*{$Z'_n$ is a polynomial of degree $n-1$ in each variable}
By applying the residue formula to Eq.~\eqref{Z'n} we obtain
\begin{multline}
Z'_n=(-1)^{\binom{n}{2}}\sum_{\substack{K=(k_1,\ldots,k_n)\\{k_l\neq k_m\ \textrm{if}\ l\neq m}\\{k_l\leq2l-1}}} (-1)^{s(K)}\prod_{l<m}(qz_{k_l}-q^{-1}z_{k_m})\label{Z'p}\\
\times\frac{\displaystyle\prod_{\substack{i<j\\i\notin K\textrm{\ or\ }(i=k_l \textrm{\ and\ } j< 2l-1)}}(qz_i-q^{-1}z_j)\prod_{2i-1\ne k_i}(qz_{2i-1}-q^{-1}z_{k_i})}{\displaystyle \prod_{\substack{i\leq2j-1\\i\notin K\textrm{\ or\ } i>k_j}}(z_{k_j}-z_i)}
\end{multline}
where $(-1)^{s(K)}$ is the sign of the permutation that orders the $k_i$.
It is enough to prove that $\lim_{z_{k_j}\rightarrow z_i}Z_n'$ exists; the
verification is a tedious but easy calculation (see \cite{DFZJ-qKZ-TSSCPP}
for a similar check).

We can now consider the leading term in each variable $z_i$
in the summation of Eq.~\eqref{Z'p}, depending on whether $i\in K$ or not;
in both cases we find a degree $n-1$.

\subsubsection*{$Z'_n$ satisfies the wheel condition}

Using the formula~\eqref{Z'p}, we can verify that $Z'_n$ is zero at $z_k=q^2z_j=q^4z_i$ for all $k>j>i$:
In fact, the term $\prod_{s<r\textrm{\ and\ } s\notin K}(qz_s-q^{-1}z_r)$ implies that $i$ and $j\in K$. As a consequence of the term  $\prod_{l<m}(qz_{k_l}-q^{-1}z_{k_m})$, we must have $i=k_m$ and $j=k_l$ with $l<m$, but, in this case, $j\leq2l-1<2m-1$ proving that $Z'_n$ satisfies the ``wheel condition''.

\subsubsection*{Recursion relation}
We show that $Z'_n$,
at $q=e^{2\pi i/3}$, satisfies a weaker form of recursion relation~\eqref{strogrecur}. Let $j$ be an integer between $1$ and $2n-1$
and evaluate $Z'_n$ at $z_{j+1}=q^2 z_j$.
We will perform the calculation for $j$ even.

If we look at formula~\eqref{Z'p} it is straightforward that all terms are zero except for $j=k_m$ and $j+1\ge2m-1$, i.e. $j=k_m=2m-2$. Using the fact that $z_{j+1}=q^2z_j$, we can derive
\begin{align*}
Z'_{n|z_{j+1}=q^2z_j}=&\prod_{i<j}(qz_i-q^{-1}z_j)(qz_i-qz_j)\prod_{k>j+1}(qz_j-q^{-1}z_k)(z_j-q^{-1}z_k)(-1)^{\binom{n}{2}}\\
&\times\prod_{i<k\neq j, j+1}(qz_i-q^{-1}z_k)\oint\ldots\oint\prod_{l}\frac{dw_l}{2\pi i} \prod_{l\neq m}(qz_{2l-1}-q^{-1}w_l)\\
&\times\frac{\prod_{l<p\neq m}(w_p-w_l)(qw_l-q^{-1}w_p)(z_j-q^{-1}z_j)}{\prod_{l\neq m}\prod_{\substack{i\leq 2l-1\\i\neq j,j+1}} (w_l-z_i) \prod_{\substack{i\geq 2l-1\\i\neq j,j+1}} (qw_l-q^{-1}z_i)}\\
&\times\prod_{n>m}\frac{(w_n-z_j)(qz_j-q^{-1}w_n)}{(w_n-z_j)(w_n-q^2z_j)}\prod_{l<m}\frac{(z_j-w_l)(qw_l-q^{-1}z_j)}{(qw_l-q^{-1}z_j)(qw_l-qz_j)}\\
&\times\frac{1}{(z_j-q^2z_j)\prod_{i<j}(z_j-z_i)\prod_{k>j+1}(qz_j-q^{-1}z_k)}
\end{align*}
After multiple cancellations we get:
\begin{equation}
Z'_n(\ldots,z_j,z_{j+1}=q^2z_j,\ldots)=\prod_{i\neq j,j+1}(q z_j-z_i)Z'_{n-1}
(z_1,\ldots,z_{j-1},z_{j+2},\ldots,z_{2n})
\end{equation}
The formula actually holds for both parities of $j$; the proof for $j$ odd is similar.

\subsubsection*{Calculating $Z'_n$ at the $c_n$ points}
Using the formula above, we can easily calculate $Z'_n$ at the $c_n$ points of the lemma. One can always choose 
two consecutive variables which are $(q^{-1},q)$ and apply the recursion relation above:
\begin{align*}
Z'_n(\ldots,z_j=q^{-1},z_{j+1}=q^2z_j=q,\ldots)&=\prod_{i\neq j,j+1}(1-z_i)Z'_{n-1}\\
&=(1-q)^{n-1}(1-q^{-1})^{n-1}Z_{n-1}'
\end{align*}
The second equality uses the fact that there is the same number of $\epsilon_i=1$ and $\epsilon_i=-1$. Since we have $Z'_1=1$, we obtain:
\[Z'_n=3^{\binom{n}{2}}\]

We finally conclude, by applying lemma \ref{wheel}, that
\[
Z_n=Z'_n
\]

Starting from our new integral formula for the
partition function of the 6-Vertex model~\eqref{Z'n}, we are now in a position
to calculate
\[\frac{(q^{2}(q+x)(1+qy))^{n-1}}{(q-q^{-1})^{n(n-1)}}Z_{n}(\frac{1+qx}{q+x},1,\ldots,1,\frac{q+y}{1+qy})\]
After some tedious computations and using new variables
\[
 u_i  =  \frac{w_i - 1}{q w_i -q^{-1}}
\]
we obtain:
\[(y+x-yx)\oint\ldots\oint\prod_{l}^{n}\frac{du_{l}}{2\pi i}\frac{1}{u_{l}^{2l-2}}\frac{\prod_{l<m}(u_{m}-u_{l})(1+u_{m}+u_{m}u_{l})}{(1+u_{l}-x)(1+u_{l}(1-y))}\prod_{j=2}^{n}\left(1+u_{j}\right)\]
where the integral contours surround counterclockwise 
$u_i=0$ and $u_i=x-1$ (and not $1/(y-1)$).

To simplify our calculation we integrate on $u_1$:
\begin{equation}
\tilde{A}_{n}(x,y)=\oint\ldots\oint\prod_{l=2}^{n}\frac{du_{l}}{2\pi i}\frac{(1+u_{l})(1+xu_{l})}{u_{l}^{2l-2}(1+u_{l}(1-y))}\prod_{l<m}^{n}(u_{m}-u_{l})(1+u_{m}+u_{m}u_{l}) \label{A}
\end{equation}
where the contours surround the remaining poles at $u_i=0$ only.

\subsection{Integral formula for refined NILP counting}
We shall derive a contour integral formula for the generating polynomial ${N}_{10}'(t_0,t_1,\ldots,t_{n-1})$ of our NILPs with a weight $t_i$ per vertical step in the $i^{\rm th}$ slice (between $y=1-i$ and $y=-i$). We use the Lindstr\"{o}m--Gessel--Viennot formula~\cite{Lind,GV} (see also the third chapter of \cite{Bressoud}):
\begin{equation}
N_{10}'(t_0,t_1,\ldots,t_{n-1})=\sum_{\substack{{1=r_1 < \ldots < r_{n-1}}\\{r_i\le 2i+1}\\r_{i+1}-r_i\ \textrm{odd}}} \det [\mathcal{P}_{i,r_j}]
\end{equation}
where $\mathcal{P}_{i,r}$ is the weighted sum over all possible lattice paths from $(i,-i)$ to $(r+1,1)$. Such paths counts with $r-i+1$ diagonal steps and $2i-r$ vertical ones, hence:
\begin{equation}
\mathcal{P}_{i,r}=\sum_{0\leq i_1<\ldots<i_{2i-r}\leq i} \prod _{l=1} ^{2i-r} t_{i_l} = \prod _{k=0} ^{i} \left. (1+t_k u)\right|_{u^{2i-r}}
\end{equation}
where the subscript $u^{2i-r}$ stands for the coefficient of the corresponding power of $u$ in the polynomial. 

We can reintroduce the path beginning at $(0,0)$ and rewrite the equation as a contour integral:
\[
{N}_{10}'(t_0,t_1,\ldots,t_{n-1})= \oint\ldots\oint \prod _{i=1} ^n\frac {du_i}{2\pi i u_i^{2i-1}}\prod _{k=0} ^{i-1} (1+t_k u_i) \sum_{\substack{{0=r_0 < r_1 < \ldots < r_{n-1}}\\r_{i+1}-r_i\ \textrm{odd}}} \det [u_i^{r_{j-1}}]
\]
where the paths of integrations are small counterclockwise circles around zero.

The last sum can be evaluated as a standard result for the sum over all Schur functions corresponding to \emph{even} partitions (see exercise 4.3.9 in \cite{Bressoud}):
\begin{equation}
\sum_{\substack{{0=r_0 < r_1 < \ldots < r_{n-1}}\\r_{i+1}-r_i\ \textrm{odd}}} \det [u_i^{r_{j-1}}]=\frac{\prod _{j>i} (u_j-u_i)}{\prod _{j\ge i} (1-u_j u_i)}
\end{equation}
where we have relaxed the condition $r_0=0$ into $r_0\ge 0$ and even, since this does not affect the integral.

The integral can thus be transformed as follows:
\begin{equation}
{N}_{10}'(t_0,t_1,\ldots,t_{n-1})=\oint\ldots\oint\prod_{i=1}^n\frac{du_i}{2\pi i u_i^{2i-1}}\frac{1}{1-u^2_i}\prod_{k=0}^{i-1}(1+t_k u_i)\prod_{j>i}\frac{u_j-u_i}{1-u_j u_i}
\end{equation}

We are mainly interested in the case where $t_0=t$, $t_1=s$ and all the others $t_i$ equal $1$. In this case, we rewrite the equation:
\begin{multline}
{N}_{10}'(t,s,1,\ldots,1):=U^{0,1}_n(t,s) \label{U}
\\=\oint\ldots\oint\prod_{i=1}^n\frac{du_i}{2\pi i u_i^{2i-1}}\frac{1}{1-u^2_i}(1+t u_i)(1+s u_i)_{\hat{1}}(1+u_i)^{i-2}_{\hat{1}}
\prod_{j>i}\frac{u_j-u_i}{1-u_j u_i}
\end{multline}
where $\hat {1}$ means that we exclude the term corresponding to $u_1$. 

\subsection{Equality of integral formulae}

At this point, we have two integral expressions, $A_n(x,y)$ (in equation~\eqref{A}) and $U^{0,1}_n(x,y)$ (in equation~\eqref{U}) and we want to prove that they are the same. The first step is to integrate over $u_1$ the expression~\eqref{U}:
\begin{equation}
 U^{0,1}_n(x,y)=\oint\ldots\oint\prod_{i=2}^{n}\frac{du_{i}}{2\pi iu_{i}^{2i-2}}(1+xu_{i})(1+yu_{i})(1+u_{i})^{i-2}\frac{\prod_{i<j}(u_{j}-u_{i})}{\prod_{i\le j}(1-u_{j}u_{i})} \label{U-1}
\end{equation}

At this stage we use the following identity:
\begin{multline}
\oint\ldots\oint\prod_{i}\frac{du_{i}}{2\pi i}\frac{\varphi(u)}{u_{i}^{2i}}\prod_{i<j}(u_{j}-u_{i})(1+\tau u_{j}+u_{i}u_{j})\\=\oint\ldots\oint\prod_{i}\frac{du_{i}}{2\pi i}\varphi(u)\frac{(1+\tau u_{i})^{i-1}}{u_{i}^{2i}}\frac{\prod_{i<j}(u_{j}-u_{i})}{\prod_{i\le j}(1-u_{i}u_{j})}
\label{zeilid}
\end{multline}
for any $\varphi(u)$ completely symmetric in $(u_1,u_2,\ldots,u_n)$ and without poles in a neighborhood of zero.
This was conjectured in~\cite{DFZJ-qKZ-TSSCPP} and proved in
\cite{Zeil}. We present in appendix~\ref{Zeq} an independent proof of a stronger
formula that implies Eq.~(\ref{zeilid}).

If we shift the indexes $(i-1)\rightarrow i$, consider $\tau=1$ and set $\varphi(u)=\prod_{i=1}^{n-1}(1+xu_i)(1+yu_i)$ we can apply the equality:
\begin{equation}
U^{0,1}_n(x,y)=\oint\ldots\oint\prod_{i=2}^{n}\frac{du_{i}}{2\pi  iu_{i}^{2i-2}}(1+xu_{i})(1+yu_{i})\prod_{i<j}(u_{j}-u_{i})(1+u_{j}+u_{j}u_{i})
\end{equation}

Now we remark that the two integrals are the same, except for
the pieces
$\frac{(1+u_{l})}{(1+u_{l}(1-y))}$ versus $1+yu_{l}$. 
Unsurprisingly, we find that is possible to write both integrals as special cases of the same integral:
\begin{equation}
 I_n(x,y)=\oint\ldots\oint\prod_{l=1}^{n-1}\frac{du_{l}}{2\pi i}\frac{(1+u_{l}+a_{l}u_{l}^{2})(1+xu_{l})}{u_{l}^{2l}(1+u_{l}(1-y))}\prod_{l<m}^{n-1}(u_{m}-u_{l})(1+u_{m}+u_{m}u_{l}) \label{I}
\end{equation}
which takes the value of $\tilde{A}_{n}(x,y)$ if $a_{l}=0$ for all $l$ and takes the value of $U^{0,1}(x,y)$
if $a_{l}=y(1-y)$ for all $l$.

More surprising is the fact that $I_n$ does not depend on the $a_i$. 
We shall show by induction on $i$ that $I_n$ is independent of $a_i$,
noting that it is a polynomial in $a_i$ of degree at most 1.

Let us first differentiate $I_n$ with respect to $a_1$:
 \begin{align*}\frac{d}{da_{1}}I_n(x,y) =&\oint\frac{du_{1}}{2\pi i}\frac{(1+xu_{1})}{(1+u_{1}(1-y))}\\&\times\oint\ldots\oint\prod_{l=2}^{n-1}\frac{du_{l}}{2\pi i}\frac{(1+u_{l}+a_{l}u_{l}^{2})(1+xu_{l})}{u_{l}^{2l}(1+u_{l}(1-y))}\prod_{m<l}^{n-1}(u_{l}-u_{m})(1+u_{l}+u_{m}u_{l})
\end{align*}
but, this integral has no poles at $u_1$ so it vanishes. 

Let us now assume by induction hypothesis that $I_n$ does not depend on the first $(i-1)$ $a_j$, and prove that the expression~\eqref{I} does not depend on $a_i$ either. As the integral does not depend on $a_j$ for all $j<i$ we can set all $a_j=0$ (for $j<i$). 

If we differentiate now with respect to $a_i$ and look at what happens in the integration up to $u_i$. We find an expression of the type:
\begin{equation}
 J_{i}  =  \oint\frac{du_{i}}{2\pi iu_{i}^{2i-2}}\oint\cdots\oint\prod_{j=1}^{i-1}\frac{du_{j}}{2\pi i}\frac{1+u_{j}}{u_{j}^{2j}}\Theta_{i}A_{i}
\end{equation}
where $A_i$ is some anti-symmetric function in the $u_j$ for all $j\le i$ without any poles in the integration domain, and $\Theta_{i}=\prod_{j<i}(1+u_{i}+u_{j}u_{i})$.

To prove that this integral is always zero we shall proceed once again by induction. The first one, $J_1$, is zero because it has no poles:
\begin{equation}
 J_1=\oint \frac{du_1}{2\pi i}A_1(x_1)=0
\end{equation}

Let $J_{i-1}=0$. All the poles are at $0$, the $A_i$ is anti-symmetric between $u_i$ and $u_{i-1}$, so we can take advantage of the fact that the $u_i$ appears with the same degree as $u_{i-1}$ in the denominator to erase all the symmetric terms in the expression $(1+u_{i-1})(1+u_i+u_{i-1} u_1)$ and get $u_i u^2_{i-1}$:
\begin{equation}
 J_{i} = \oint\frac{du_{i}}{2\pi iu_{i}^{2i-3}}\oint\frac{du_{i-1}}{2\pi iu_{i-1}^{2i-4}}\oint\ldots\oint\prod_{j=1}^{i-2}\frac{du_{j}}{2\pi i}\frac{1+u_{j}}{u_{j}^{2j}}\hat{\Theta}_{i}A_{i} \label{J_s}
\end{equation}
where the hat in $\hat{\Theta}_{i}$ means that the term $(1+u_{i}+u_{i-1}u_{i})$ is skipped\footnote{Note that $\hat{\Theta}_i$ is symmetric between $u_{i-1}$ and $u_i$.}. The integral does not have yet the desired form, i.e.\ $J_{i-1}$, it is missing the term $(1+u_{i}+u_{i-1}u_{i})$, so we add and subtract it:
\begin{align*}
J_{i}  =& \oint\ldots\oint\frac{du_{i}}{2\pi iu_{i}^{2i-3}}\frac{du_{i-1}}{2\pi iu_{i-1}^{2i-4}}\prod_{j=1}^{i-2}\frac{du_{j}}{2\pi i}\frac{1+u_{j}}{u_{j}^{2j}}(1+u_{i}+u_{i}u_{i-1}-u_{i}-u_{i}u_{i-1})\hat{\Theta}_{i}A_{i}\\
 =& \oint\ldots\oint\frac{du_{i}}{2\pi iu_{i}^{2i-3}}\frac{du_{i-1}}{2\pi iu_{i-1}^{2i-4}}\prod_{j=1}^{i-2}\frac{du_{j}}{2\pi i}\frac{1+u_{j}}{u_{j}^{2j}}\Theta_{i}A_{i}\\
& - \oint\ldots\oint\frac{du_{i}}{2\pi iu_{i}^{2i-4}}\frac{du_{i-1}}{2\pi iu_{i-1}^{2i-4}}\prod_{j=1}^{i-2}\frac{du_{j}}{2\pi i}\frac{1+u_{j}}{u_{j}^{2j}}(1+u_{i-1})\hat{\Theta}_{i}A_{i}
\end{align*}
The first term is already in the form of $J_{i-1}$. The second term is almost symmetric between $u_{i}$ and $u_{i-1}$, using the same method as in~\eqref{J_s} we can transform $-(1+u_{i-1})$ to $1+u_{i}+u_{i}u_{i-1}$; in this way, we recover the symmetry needed so that we can write $J_i$ as an integral in $u_i$ of some function multiplied by $J_{i-1}$, which is zero. As a consequence $J_i$ is also zero for all $i$, i.e. $I_n(x,y)$ does not depend on any $a_i$. We conclude
that
\[\tilde{A}_n(x,y)=U^{0,1}_n(x,y)\]
\begin{flushright}
 \qedsymbol
\end{flushright}

\appendix

\section{Formulating the conjecture directly in terms of TSSCPPs} 
We have used the NILP formulation throughout this paper (in particular, to prove the main theorem), whereas Mills, Rumsey and Robbins use the language of TSSCPPs. In \ref{U-u} we first describe the theorem in a more general form, and then prove that we can reduce it to the one presented in~\ref{conj}. 
We then reformulate in~\ref{Mills2} 
our theorem in the language of~\cite{MRR-TSSCPP,rob-story}.

\subsection{Extending the theorem} \label{U-u}

Let $A_n(x,y)$ and $\tilde{A}_n(x,y)$ be the same as defined in~\ref{A-n}. We use the same NILPs with the extra-step as in~\ref{N-X}. 

\begin{figure}
\centering 
\psfrag{a}[0][0][1][0]{$u_6^0$}\psfrag{b}[0][0][1][0]{$u_6^3$}
\includegraphics[scale=0.7]{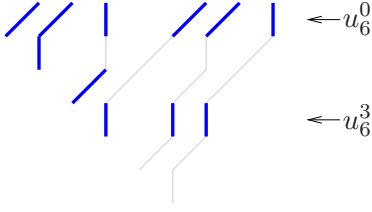}
\caption{Let $\alpha$ be the NILP represented here. In order to calculate $u_6^0(\alpha)$ and $u_6^3(\alpha)$ we highlight the extra-steps and the $\max\{1,t-3+1\}$-th step of the path starting at $(t,-t)$.\label{u_0_3} Here we have $u_6^0(\alpha)=2$ and $u_6^3(\alpha)=4$.}
\end{figure}

We now introduce a function $u^k_n(\alpha)$, where $\alpha$ is a NILP, which counts the number of vertical steps in the extra-step if $k=0$; otherwise it counts the number of vertical steps in the $\max\{1,t-k+1\}$-th step of the path starting at $(t,-t)$, as shown on figure~\ref{u_0_3}.

We can next define the function $U^i_n(x)$:
\begin{equation}
 U^i_n(x):= \sum_k U^i_{n,k} x^k := \sum_{\alpha} x^{u_n^i(\alpha)}
\end{equation}
and more complex functions $U^{i,j}_n(x,y)$:
\begin{equation}
 U^{i,j}_n(x,y):= \sum_k U^{i,j}_{n,k,l} x^k y^l := \sum_{\alpha} x^{u_n^i(\alpha)} y^{u_n^j(\alpha)}
\end{equation}

We could generalize these even more, introducing more indices, but this is general enough for our purposes. With these new functions we can rewrite our theorem:

\begin{conj*}
 \begin{align}
  \tilde{A}_{n}(x,y)&=U_{n}^{0,j}(x,y) \label{T1}\\A_{n}(x,y)&=U_{n}^{1,i}(x,y) \label{T2}
 \end{align}
where $j={1,2\ldots}$ and $i={2,3\ldots}$. If we choose $U_{n}^{0,1}$ we have the theorem as stated before.
\end{conj*}

On order to reduce this 
to our previous result, it is enough to prove that $U_n^{0,i}$ does not depend in $i$ and that $U^{0,i}_{n,k,j}=U^{1,i}_{n,n-k-1,j}$ (for $i\ge 2$). 

\subsubsection*{$i$ Independence of $U_n^{0,i}$}

\begin{figure}

\centering\includegraphics[scale=0.3]{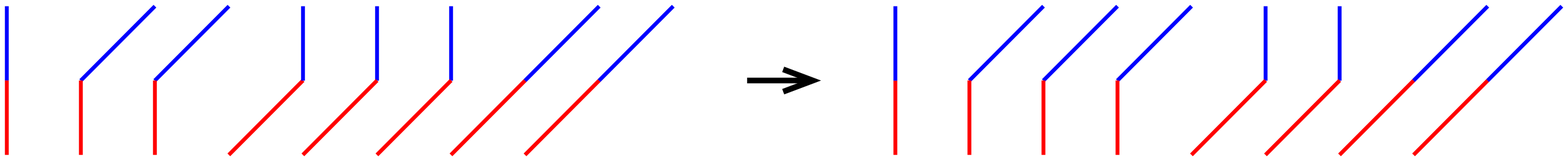}
\caption{We can group the double steps in islands, such that all the starting points (of the double steps) are consecutive. These doubles steps are, necessarily, ordered in $r$ double vertical steps, $s$ vertical-diagonal steps, $t$ diagonal-vertical steps and $u$ double diagonal steps. Our function $g$ interchanges $s$ with $t$ at each island, so that we interchange the number of vertical steps between the two rows.\label{Usym3}}

\end{figure}

For the first equality we introduce a function $g$ as explained on figure~\ref{Usym3}. This function interchanges the number of vertical steps in two consecutive rows leaving invariant all the other rows. This function has the important property $g\circ g=Id$. So, it is straightforward from this that $U_n^{0,i}=U_n^{0,i+1}$, with $i$ greater than $0$.

\subsubsection*{$U_{n,k,j}^{0,i}=U_{n,n-k-1,j}^{1,i}$ for $i>1$}
\begin{figure}

\centering\includegraphics[scale=0.3]{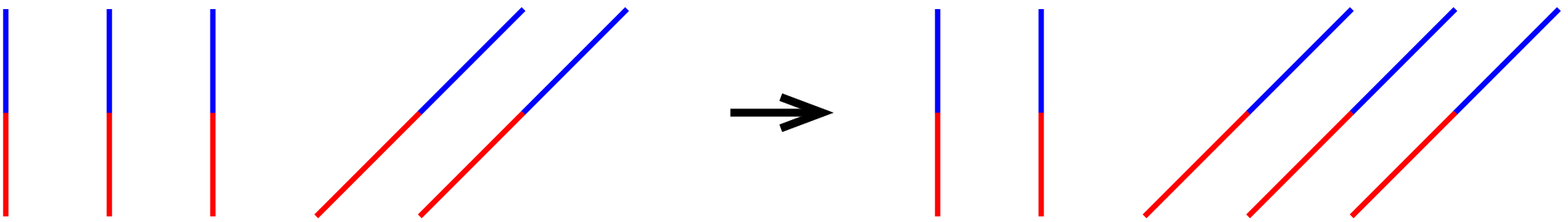}

\caption{In order to satisfy the extra-step rules we can only build two type of islands, one made of $r$ double vertical steps and $s$ double diagonal steps, and the other type made of $t$ vertical-diagonal steps and $u$ diagonal-vertical steps. Our function $h$ interchange simply $r$ with $s$. It is important to note that the first path is always invariant under $h$ (it is always of the type vertical-diagonal or the inverse).\label{Usym2}}

\end{figure}

The proof follows the same structure as the former. We construct again a function $h$ such that $h\circ h=Id$, which interchanges the number of vertical steps at the extra-step with the number of diagonal steps at the last step (before the extra-step). This function is obviously a bijection and it leaves invariant all the rows except the last one and the extra one because it is applied at the top of the diagrams as can be seen on figure~\ref{Usym2}. An important remark is that the first path is always invariant under $h$ because it is of the type vertical-diagonal or diagonal-vertical. This proves our equality.

In conclusion, all these variations (\eqref{T1} and \eqref{T2}) are truly the same, and we can concentrate on only one version.

\subsection{The conjecture in terms of TSSCPPs} \label{Mills2}

Mills, Robbins and Rumsey conjectured this theorem by means of TSSCPPs, not NILPs, but behind the different formulations lies the same result. To show that, we describe some of the content of~\cite{MRR-TSSCPP} and explain the equivalence.

\begin{figure}
\centering\includegraphics[scale=0.8]{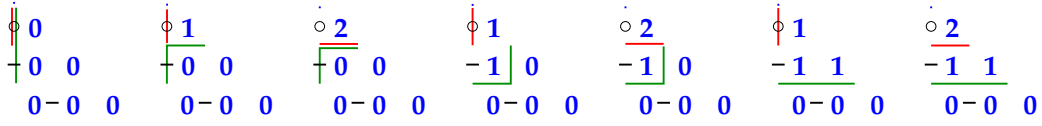}\caption{We can see on this figure what the function $u^2_3$ counts. The signs minus represents the part: $a_{t,t-k}-a_{t,t-k+1}$, so they count the vertical steps, and the little circles represents $\#\{a_{t,n+1}\ |\ a_{t,n+1}<2n-t\}$. If we stretch our diagrams to obtain the NILPs we recover our definition of $u_n^k$.\label{NILPX_R4}}

\end{figure}

Recall that TSSCPPs can be represented as $2n\times 2n$ matrices $a$,
as in Eq.~\eqref{TSSCPP}.
In \cite{MRR-TSSCPP} is 
introduced a quantity which we shall denote by $u_n^k (a)$,
and which depends on the upper-left $n\times n$ submatrix of $a$:
\begin{equation}
 u_{n}^{k}(a)=\sum_{t=1}^{n-k+1}(a_{t,t+k-1}-a_{t,t+k})+\sum_{t=n-k+2}^{n}\#\{a_{t,n}\ |\ a_{t,n}>2n-t+1\}
\end{equation}
where $\#$ means cardinality, and where conventionally,
$a_{t,n+1}:=2n-t+1$ in this equation.
Also defined is the function:
\begin{equation}
 U_{n}^{i,\ldots k}(x,\ldots,z)=\sum_{a}x^{u_{n}^{i}(a)}\ldots z^{u_{n}^{k}(a)}\qquad\ \textrm{for\ all}\ i,\ldots,k\in\{1,\ldots,n+1\}
\end{equation}
We claim that these are our functions $u$ and $U$ defined above. To make the connection, reexpress this function in terms of the lower-right
$n\times n$ submatrix of $a$:
\begin{equation}
 u_{n}^{k}(a)=\sum_{t=n+k}^{2n}(a_{t,t-k}-a_{t,t-k+1})+\sum_{t=n+1}^{n+k-1}\#\{a_{t,n+1}\ |\ a_{t,n+1}<2n-t\}
\end{equation}
where we replace $a_{t,n}$ with $2n-t$. What this function counts is described on figure~\ref{NILPX_R4}. Finally, if we shift the diagrams to obtain NILPs we recover our functions $U_n^k$ as expected.

As a final remark, in the article~\cite{rob-story} three functions are defined: $f_1$, $f_2$ and $f_3$ and the conjecture is stated with any two of them. In fact, $f_1$ is connected with the $u_n^0$, $f_2$ with the $u_n^1$ and $f_3$ with $u_n^n$, as can be seen using the same procedure.

\section{Properties of the 6-Vertex model partition function} \label{Z=Z' Z}
Let us consider, as in section~\ref{6-V pf},
the 6-Vertex model with Domain Wall Boundary Conditions.
Let $\tilde Z_n$ be its partition function
(with Boltzmann weights given by Fig.~\ref{poids}), and $Z_n$ to
be $\tilde Z_n$ divided by the normalization factor
$(-1)^{n(n-1)/2}{(q^{-1} - q)}^n \prod _{i=1} ^{2n} z_i^{1/2}$.

\begin{figure}
\centering
\psfrag{a}[0][0][1][0]{$z_1$}\psfrag{c}[0][0][1][0]{$z_2$}\psfrag{b}[0][0][1][0]{$w$}
\includegraphics[scale=1]{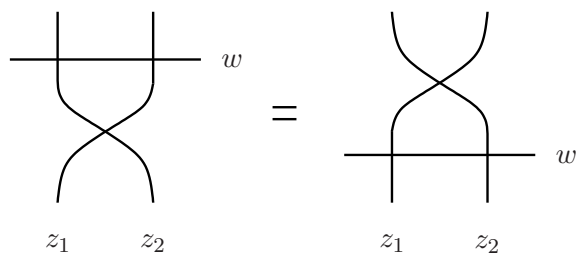}
\caption{\label{ybe}
Yang--Baxter equation. Summation over arrows of the {\it internal}\ 
edges is implied, while the external arrows are fixed and the equality holds for any choice of them.}
\end{figure}

The model thus defined satisfies the following essential property
(Yang--Baxter equation) shown on Fig.~\ref{ybe}. 
The vertex with diagonal edges is assigned
weights (the so-called $R$ matrix) which are those
of Fig.~\ref{poids} in which we have rotated the picture 45 degrees clockwise, and
with parameters $z_1$, $q^{1/2}z_2$.
parameter. In fact here we do not need
the explicit expression of the $R$ matrix, only 
that it is invariant by reversal of all arrows and
that it satisfies the ice rule i.e.\ there are as many outgoing arrows
as incoming arrows. Since the Yang--Baxter equation is invariant by change of normalization
of $R$, we can divide all weights by $b$ in such a way
that $R_{\uparrow\uparrow}^{\uparrow\uparrow}=
R_{\downarrow\downarrow}^{\downarrow\downarrow}=1$, with obvious notations.

\subsection{Korepin recursion relation}
In this paragraph, $q$ is kept arbitrary.
We shall now list the following four properties which
determine entirely $Z_n$ and only sketch
their proof (since they 
have been reproved many times since their original appearance \cite{kor}, 
see for example \cite{kup-ASM,korzj})

\begin{itemize}
\item $Z_1=1$.

This is by definition.

\item $Z_n$ is a symmetric function of the
sets of variables
$\{ z_1,\ldots,z_n\}$ and $\{ z_{n+1},\ldots,z_{2n}\}$.

It is sufficient to prove that exchange of $z_i$ and $z_{i+1}$
(for $1\le i<n$) leaves the partition function unchanged.
This can be obtained by repeated use of the Yang--Baxter property.
Multiplying the partition function by $R(z_{i+1}/z_i)$ and
noting that it is unchanged, we find
\newsavebox{\tempbox}
\newlength{\templength}
\newcommand{\vcenterbox}[1]
   {\sbox{\tempbox}{#1}
    \settowidth{\templength}{\usebox{\tempbox}}
    \parbox{\templength}{\usebox{\tempbox}}}
\begin{align*}
\tilde Z_n(\ldots, z_i,
z_{i+1},\ldots)
&=\vcenterbox{\psfrag{a}[0][0][.7][0]{$z_i$}\psfrag{c}[0][0][.7][0]{$z_{i+1}$}\includegraphics[width=3cm]{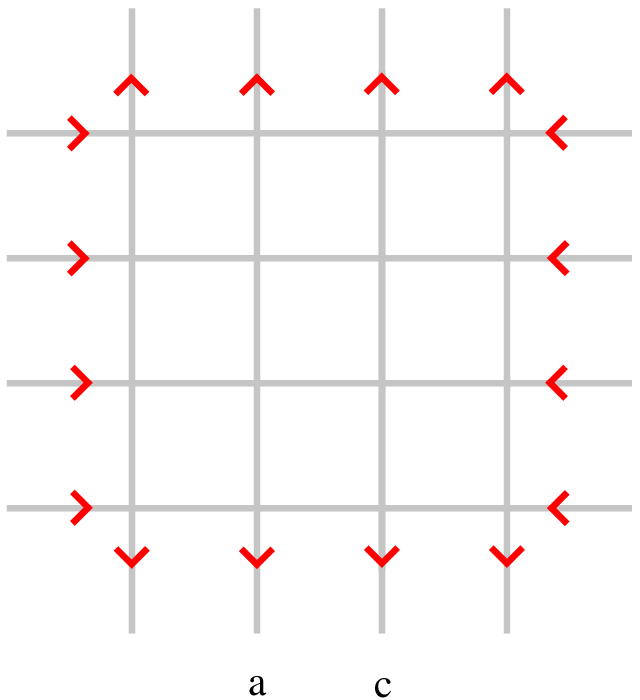}}
=\vcenterbox{\psfrag{a}[0][0][.7][0]{$z_i$}\psfrag{c}[0][0][.7][0]{$z_{i+1}$}\includegraphics[width=3cm]{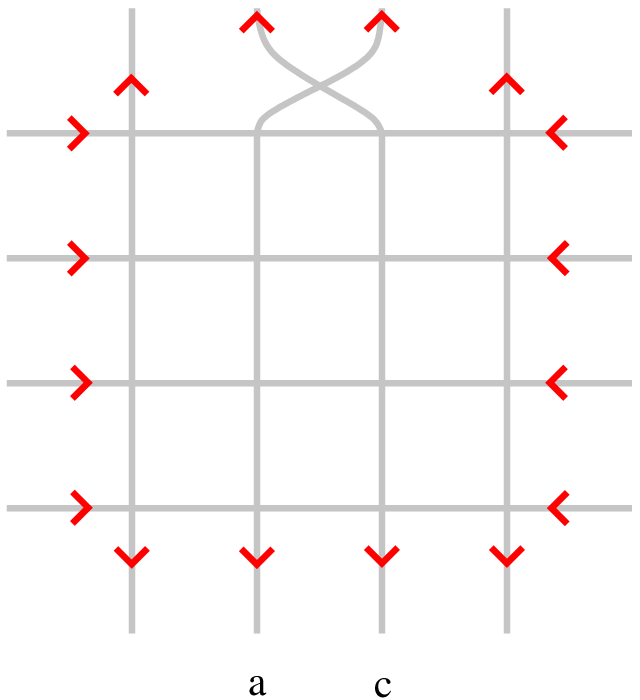}}
=\vcenterbox{\psfrag{a}[0][0][.7][0]{$z_i$}\psfrag{c}[0][0][.7][0]{$z_{i+1}$}\includegraphics[width=3cm]{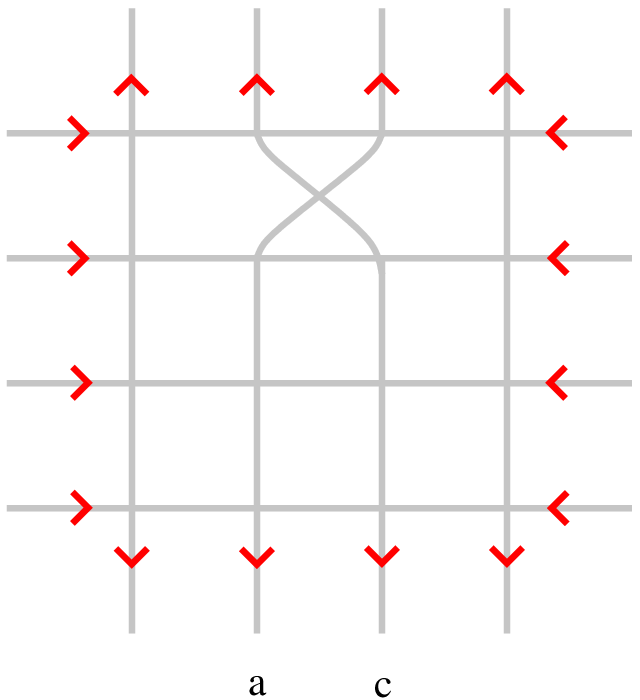}}
\\
=\cdots
&=\vcenterbox{\psfrag{a}[0][0][.7][0]{$z_i$}\psfrag{c}[0][0][.7][0]{$z_{i+1}$}\includegraphics[width=3cm]{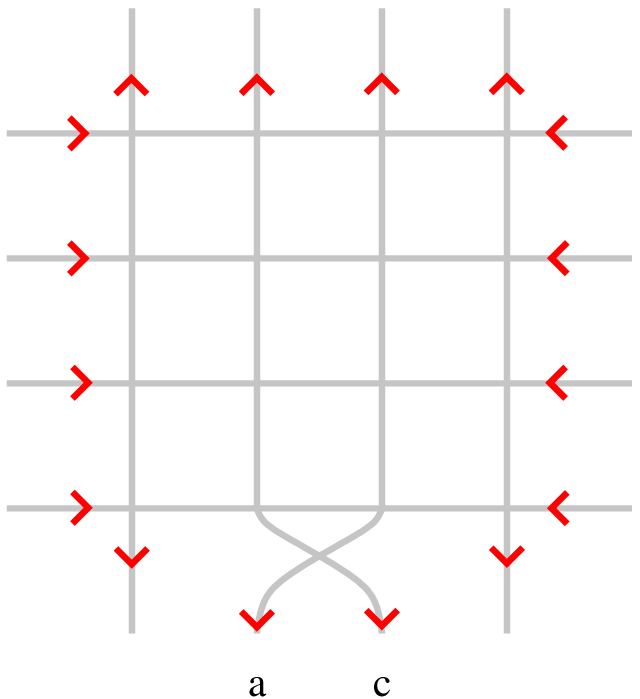}}
=\vcenterbox{\psfrag{a}[0][0][.7][0]{$z_i$}\psfrag{c}[0][0][.7][0]{$z_{i+1}$}\includegraphics[width=3cm]{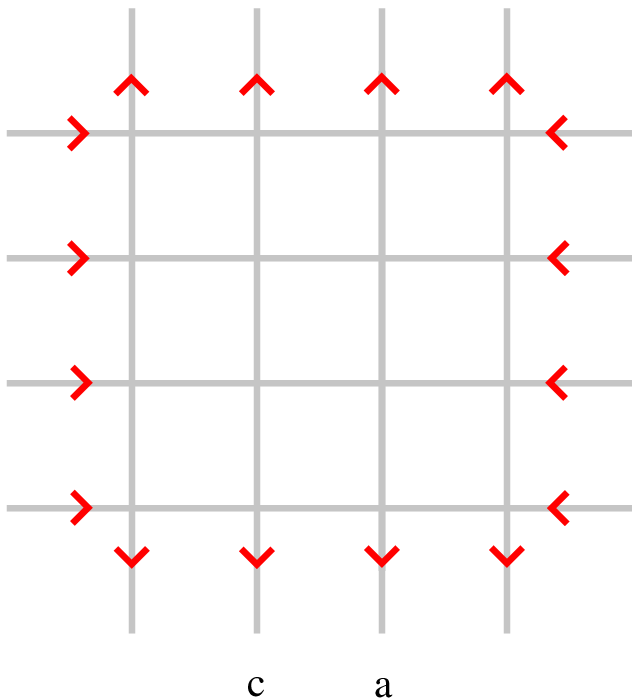}}
=\tilde Z_n(\ldots,z_{i+1},
z_i,\ldots)
\end{align*}
and similarly for the $\{ z_{n+1},\ldots,z_{2n}\}$.

\item $Z_n(z_1,\ldots,z_{2n})$ is a polynomial of degree (at most) $n-1$ in each variable.

Let us choose one configuration. Then the only weights
which depend on $z_i$ are the $n$ weights on row $i$.
Since the outgoing arrows are in opposite directions, the number of vertices
of type $c$ on this row is odd, and in particular is at least 1.
Power counting then shows that the contribution to the partition function
of any configuration is of the form $z_i^{1/2}$ times a polynomial
of $z_i$ of degree at most $n-1$. Summing over all configurations
and removing $z_i^{1/2}$ by definition of $Z_n$, we obtain the desired
property.

\item The $Z_n$ obey the following recursion relation:
\begin{multline}\label{koreprecur}
Z_{n}(z_{1},\ldots,z_n;z_{n+1}=q^{-1}z_{1},\ldots,z_{2n})\\=
q^{-n+1}\prod_{j=2}^n(z_{1}-q^{2}z_{j})
\prod_{j=n+2}^{2n}(z_1-q^{-1} z_{j})
Z_{n-1}(z_{2},\ldots,z_n;z_{n+2},\ldots,z_{2n})
\end{multline}
Since $z_{n+1}=q^{-1} z_1$ implies $a(z_{n+1},z_1)=0$,
by inspection all configurations with non-zero weights are of the
form shown on Fig.~\ref{recurs}. 
This produces the following identity for unnormalized partition functions
\begin{multline*}
\tilde Z_{n}(z_{1},\ldots,z_n;z_{n+1}=q^{-1}z_{1},\ldots,z_{2n})=
(q^{-1}-q)z_1 q^{-1/2} \\
\times \prod_{j=2}^n(q^{-3/2} z_{1}-q^{1/2}z_{j})
\prod_{j=n+2}^{2n}(q^{-1/2} z_{j}-q^{1/2}z_{1})
\tilde Z_{n-1}(z_{2},\ldots,z_n;z_{n+2},\ldots,z_{2n})
\end{multline*}
which in turns leads to the recursion relation above for the $Z_n$.

\begin{figure}
\centering
\psfrag{a}[0][0][1][0]{$z_{n+1}$}\psfrag{b}[0][0][1][0]{$z_{1}$}
\includegraphics[scale=0.7]{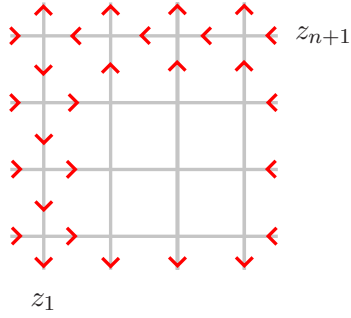}
\caption{\label{recurs}Graphical proof of the recursion relation.}
\end{figure}
\end{itemize}

Note that by the symmetry property, Eq.~\eqref{koreprecur} fixes
$Z_n$ at $n$ distinct values of $z_{n+1}=q^{-1}z_i$, $i=1,\ldots,n$.
Since $Z_n$ is of degree $n-1$ in $z_{n+1}$, it is entirely determined by it.

\subsection{Cubic root of unity case}
Let us set $q=e^{2\pi i/3}$. First, once can simplify the recursion relation 
\eqref{koreprecur} to
\[
Z_{n}(z_{1},\ldots,z_{n+1}=q^2z_{1},\ldots,z_{2n})=\prod_{j\ne 1,n+1}(q^{-2} z_{1}-z_{j})Z_{n-1}(z_{2},\ldots,z_n,z_{n+2},\ldots,z_{2n})
\]

Secondly, one wishes to show the enhanced symmetry property of $Z_n$
in the full set of variables $\{z_1,\ldots,z_{2n}\}$. For this,
it is simplest to prove Eq.~\eqref{schur}, which displays explicitly
this symmetry. Let us show that the Schur function $s_{Y_n}(z_1,\ldots,z_{2n})$
satisfies all the properties of the previous section.

$s_{Y_0}=1$ by definition.
$s_{Y_n}$ is symmetric in all variables (which is what we want to prove
for $Z_n$),
and therefore in particular symmetric in the
$\{z_1,\ldots,z_n\}$ and $\{z_{n+1},\ldots,z_{2n}\}$. 
It is a polynomial of degree $n-1$ in each variable
because the width of the Young diagram $Y_n$ is $n-1$. Finally, to obtain the recursion
relation, we note that as soon as $(z_i,z_j,z_k)=(z,q^2z,q^4z)$ for distinct $i,j,k$, the three
corresponding rows in the numerator of Eq.~\eqref{schur} are linearly dependent so that
the numerator vanishes while the denominator does not. Thus, at $z_j=q^2 z_i$, $i\ne j$,
\[
s_{Y_n}(z_{1},\ldots,z_{j}=q^2z_{i},\ldots,z_{2n})=
\prod_{k\ne i,j}(q^{-2} z_{i}-z_{k})
Z''_{n-1}(z_{2},\ldots,\hat z_i,\ldots,\hat z_j,\ldots,z_{2n})
\]
where $Z''_{n-1}$ does not depend on $z_i$ because the $2n-2$ prefactors exhaust the degree
in $z_i$.

Now set $z_i=0$: the Schur function has $2n-2$ remaining arguments, so the full column
of length $2n-2$ can be factored out and we are left with the Young diagram $Y_{n-1}$:
\[
s_{Y_n}(z_{1},\ldots,z_i=0,\ldots,z_{j}=0,\ldots,z_{2n})=
\prod_{k\ne i,j}z_{k}\ 
s_{Y_{n-1}}(z_{2},\ldots,\hat z_i,\ldots,\hat z_j,\ldots,z_{2n})
\]
By comparison, we conclude that $Z''_{n-1}=s_{Y_{n-1}}$,
so that $s_{Y_n}$ satisfies the desired recursion relation.

We conclude that $s_{Y_n}$ satisfies all the properties of the
previous section, which determine
uniquely $Z_n$. Thus, Eq.~\eqref{schur} holds.

\section{The space of polynomials satisfying the wheel condition} \label{Z=Z' lemma}

In order to prove that $Z'_n$ (defined in~\eqref{Z'n}) is the partition function of the 6-V model, we need to prove lemma~\ref{wheel}. That is, a polynomial $P$ of degree (at most) $n-1$ in each variable $z_1,\ldots,z_{2n}$ satisfying the ``wheel condition'' is entirely determined by its values at the following specializations: $(q^{\epsilon_1},\ldots,q^{\epsilon_{2n}})$ for all possible choices of $\{\epsilon_i=\pm 1\}$ such that $\sum_{i=1}^{2n} \epsilon_i=0 $ and $\sum_{i=1}^j\epsilon_i\leq 0$ for all $j\leq2n$ (these are just increments of Dyck paths).

Or equivalently, if a polynomial satisfies these conditions and is zero at all the specializations, then it is identically zero. For example, at $n=1$ the polynomial is of degree 0 i.e.\ a constant, and as it vanishes at $(z_1,z_2)=(q^{-1},q)$ it is identically zero.

We now proceed by induction.
We suppose that the lemma is true for $n<p$. Let $\phi_p$ be a polynomial of degree $(p-1)$ at each variable which is zero at all specializations. The polynomial satisfies the ``wheel condition'' at $z_{i+1}=q^2 z_i$, so we can write
\begin{equation}
 \phi_{p}(z_1,\ldots,z_{2p})_{|z_{i+1}=q^2 z_i}=\prod_{j\neq i,i+1} (qz_i-z_j) \psi_{p-1} (z_1,\ldots,z_{i-1},z_{i+2},\ldots,z_{2p}) \label{phi-wc}
\end{equation}
where $\psi_{p-1}$ is a function of degree $p-2$ in each $z_j$ (except $z_i$ and $z_{i+1}$) which still follows the ``wheel condition''. Furthermore, let $\pi_p$ be a specialization which has $(z_i,z_{i+1})=(q^{-1},q)$ and $\pi_{p-1}'$ the same specialization but without $z_i$ and $z_{i+1}$. We apply~\eqref{phi-wc}:
\begin{equation}
 \phi_{p}(\pi_p)=(1-q)^{n-1}(1-q^{-1})^{n-1} \psi_{p-1} (\pi_{p-1}')=0
\end{equation}
The mapping $\pi_p\mapsto \pi_{p-1}$ is a bijection from Dyck paths
with $(q^{-1},q)$ at locations $(i,i+1)$ to all Dyck paths.
Thus our induction hypothesis applies, and $\psi_{p-1}=0$.

Therefore, one can write:
\begin{equation}
 \phi_p=\prod_{i=1}^{2n-1}(z_{i+1}-q^2z_i)\phi_p^{(1)}
\end{equation}
where $\phi_p^{(1)}$ is a polynomial of degree $\delta_1=\delta_{2p}=p-2$ at $z_1$ and $z_{2p}$ and $\delta_i=p-3$ at all the other variables which follows a weak version of the ``wheel condition'': 
\[
\phi^{(1)}_{p|z_k=q^2z_j=q^4z_i}=0\qquad\textrm{\ for\ all\ }k\geq j+2\geq i+4
\]
This implies:
\begin{equation}
\phi^{(1)}_{p|z_{i+2}=q^2z_i}=\prod_{j\notin [i-1,i+3]}(qz_1-z_j) \psi_p^{(1,i)}
\end{equation}
By degree counting in $z_i$ we find that they are identically zero.

Now, we can write
\begin{equation}
 \phi^{(1)}_p=\prod_{i=1}^{2n-2}(z_{i+2}-q^2z_i)\phi^{(2)}_p
\end{equation}
where $\phi^{(2)}_p$ has degree $\delta_1=\delta_{2p}=p-3$, $\delta_2=\delta_{2p-1}=p-4$ and all the others $\delta_i=p-5$.

Clearly, this procedure can be repeated; at step $r$, $\phi_p^{(r)}$ has degree:
\begin{align*}
 \delta_1&=p-r-1\\
 \delta_2&=p-r-2\\
 &\vdots\\
 \delta_r&=p-2r\\
 &\vdots\\
 \delta_i&=p-2r-1\\
 &\vdots\\
 \delta_{2p}&=p-r-1
\end{align*}
We write
\[
 \phi^{(r)}_{p|z_{i+r+1}=q^2z_i}=\prod_{j\notin[i-r,i+2r+1]}(qz_i-z_j)\psi^{(r,i)}_{p}
\]
Counting the degree in $z_i$ we conclude that $\psi^{(r,i)}_{p}=0$. So we can construct $\phi^{(r+1)}_{p}$. 

When $r\geq \frac{n}{2}$ we obtain a polynomial of negative degree which implies that the polynomial is identically zero.

Remark: What this lemma shows in other words is that the vector space
of polynomials of degree at most $n-1$ in each variable
satisfying the wheel condition is of dimension at most $c_n$.
In fact it is known to be of dimension exactly $c_n$; the
standard proof involves the fact that it is an irreducible representation
of the affine Hecke algebra, see e.g.~\cite{Pasquier,Kasatani}.

\section{An antisymmetrization formula} \label {Zeq}
The goal of this section is to prove identity \eqref{zeilid}, which allows to turn an
equation of the type~\eqref{U-1} into one of the type~\eqref{A}. 
Identity \eqref{zeilid} was conjectured by Di Francesco and Zinn-Justin
in ~\cite{DFZJ-qKZ-TSSCPP} and proved by Zeilberger \cite{Zeil}. Equivalently,
it was proved that the integrand of the l.h.s.\ without the factor $\varphi(u)$, 
once antisymmetrized
and truncated to its negative degree part (the positive powers of the $u_i$
cannot contribute to the integral), reduces to the integrand of the r.h.s.\ 
without the factor $\varphi(u)$.
Here we prove in an independent way a much stronger statement. Indeed,
here we perform the {\it exact}\/ antisymmetrization of a spectral parameter
dependent generalization of the integrand.\footnote{More precisely,
the expression we antisymmetrize is the integrand before the homogeneous
limit in which spectral parameters come in pairs $\{z,1/z\}$.}

\subsection{The general case}
Let $h_{q}(x,y)=(qx-q^{-1}y)(qxy-q^{-1})$ (and, obviously, $h_{1}(x,y)=(x-y)(xy-1)$). Let us also define
\begin{align}
f(w,z) & =  \frac{1}{z(1-q^{2}w^{2})(q^{-2}-1)}\left(\frac{1}{h_{1}(w,z)}-\frac{1}{h_{q}(w,z)}\right) \label{f}\\
 & =  \frac{1}{h_{1}(w,z)h_{q}(w,z)}
\nonumber
\end{align}
The quantity of interest is
\begin{equation}
 B_{n}(w,z)=AS\left\{ \frac{\prod_{i<j}(qw_{i}-q^{-1}w_{j})}{\prod_{i\le j}h_{1}(w_{j},z_{i})\prod_{i\ge j}h_{q}(w_{j},z_{i})}\right\} \label{BZ}
\end{equation}
where $AS(\phi)(w_1,\ldots,w_n)=\sum_{\sigma\in\mathcal{S}_n}(-1)^{|\sigma|}\phi(w_{\sigma(1)},\ldots,w_{\sigma(n)})$ 

We then claim that $B_n$ can be written as:
\begin{equation}
B_{n}(w,z)=\frac{q^{\frac{n(n-1)}{2}}\mathbf{f}_{n}}{\prod_{i<j}h_{1}(z_{i},z_{j})(1-q^{2}w_{i}w_{j})} \label{Bf}
\end{equation}
where $\mathbf{f}_{n}=\det[f(w_{i},z_{j})]_{i,j\le n}$.

Again, we prove it by induction. For $n=1$, we obtain on both sides:
\[B_1=\frac{1}{h_{1}(w_{1},z_{1})h_{q}(w_{1},z_{1})}\]

Let the equality of~\eqref{BZ}~and~\eqref{Bf} hold at $n-1$. Starting from~\eqref{BZ} and pushing $z_n$ and $w_j$ out of the anti-symmetrization we can write our equation as follows:
\begin{equation*}
B_n(w,z)=\sum_{j}(-1)^{n+j}\frac{\prod_{i\neq j}(qw_{i}-q^{-1}w_{j})}{\prod_{i}h_{1}(w_{j},z_{i})h_{q}(w_{i},z_{n})}AS\left\{ \frac{\prod_{i<k}(qw_{l}-q^{-1}w_{k})}{\prod_{i\le k}h_{1}(w_{k},z_{i})\prod_{i\ge k}h_{q}(w_{k},z_{i})}\right\} _{\hat{z}_{n}\hat{w}_{j}}
\end{equation*}
where the hat over $\hat{z_n}$ and $\hat{w_j}$ means that the terms that include them are absent from the anti-symmetrization. We use the hypothesis to replace the anti-symmetrization part:
\begin{align*}
B_{n} & =  \sum_{j}(-1)^{n+j}\frac{\prod_{i\neq j}(qw_{i}-q^{-1}w_{j})}{\prod_{i}h_{1}(w_{j},z_{i})h_{q}(w_{i},z_{n})}\frac{q^{\frac{(n-1)(n-2)}{2}}\mathbf{f}_{n-1,\hat{w}_{j}\hat{z}_{n}}}{\left(\prod_{i<k}h_{1}(z_{i},z_{k})(1-q^{2}w_{i}w_{k})\right)_{\hat{w}_{j}\hat{z}_{n}}}\\
 & =  \sum_{j}(-1)^{n+j}\frac{\prod_{i\neq j}h_{q}(w_{i},w_{j})\prod_{i\neq n}h_{1}(z_{i},z_{n})}{\prod_{i}h_{1}(w_{j},z_{i})h_{q}(w_{i},z_{n})}\frac{(-1)^{n-1}q^{\frac{n(n-1)}{2}}\mathbf{f}_{n-1,\hat{w}_{j}\hat{z}_{n}}}{\left(\prod_{i<k}h_{1}(z_{i},z_{k})(1-q^{2}w_{i}w_{k})\right)}
\end{align*}

The idea now is to rewrite this expression under the form
$\sum_{j}(-1)^{n+j}\mathbf{f}_{n-1,\hat{w}_{j}\hat{z}_{n}}\sum_{i}g_{i}f(w_{j},z_{i})$
for some functions $g_i$.
Indeed, using the fact that $\mathbf{f}_n$ is a determinant, we would get
\[
\sum_{j}(-1)^{n+j}\mathbf{f}_{n-1,\hat{w}_{j}\hat{z}_{n}}\sum_{i}g_{i}f(w_{j},z_{i})
=\sum_{j}(-1)^{n+j}\mathbf{f}_{n-1,\hat{w}_{j}\hat{z}_{n}}g_{n}f(w_{j},z_{n})=\mathbf{f}_{n}g_{n}\]

One can guess the form of $g_i$:
\[g_{i}=\frac{\prod_{j\neq i,n}h_{1}(z_{j},z_{n})\prod_{j}h_{q}(w_{j},z_{i})}{\prod_{j\neq i,n}h_{1}(z_{i},z_{j})\prod_{j}h_{q}(w_{j},z_{n})}\]
One can verify this decomposition directly. Equivalently, it can be written as
\begin{equation}
\sum_{i}\frac{\prod_{j\neq i,n}h_{1}(z_{j},z_{n})\prod_{j}h_{q}(w_{j},z_{i})}{\prod_{j\neq i,n}h_{1}(z_{i},z_{j})\prod_{j}h_{q}(w_{j},z_{n})}f(w_{k},z_{i})=\frac{\prod_{i\neq k}h_{q}(w_{i},w_{k})\prod_{i\neq n}h_{1}(z_{i},z_{n})}{\prod_{i}h_{1}(w_{k},z_{i})h_{q}(w_{i},z_{n})}
\end{equation}
or, by multiply both sides with $\prod_{i}h_{1}(w_{k},z_{i})h_{q}(w_{i},z_{n})$ to obtain polynomials of $w_k$ of degree $2(n-1)$:
\begin{equation}
\sum_{i}\frac{\prod_{j\neq i,n}h_{1}(z_{j},z_{n})}{\prod_{j\neq i,n}h_{1}(z_{i},z_{j})}\prod_{j\neq k}h_{q}(w_{j},z_{i})\prod_{j\neq i}h_{1}(w_{k},z_{j})=\prod_{i\neq k}h_{q}(w_{i},w_{k})\prod_{i\neq n}h_{1}(z_{i},z_{n})
\end{equation}

It is enough to prove that this equation is the same in all points $w_k=z_i$ and $w_k=z_i^{-1}$. In the first case we have: 
\begin{align*}
\frac{\prod_{j\neq i,n}h_{1}(z_{j},z_{n})}{\prod_{j\neq i,n}h_{1}(z_{i},z_{j})}\prod_{j\neq k}h_{q}(w_{j},z_{i})\prod_{j\neq i}h_{1}(z_{i},z_{j}) & =  \prod_{j\neq k}h_{q}(w_{j},z_{i})\prod_{j\neq n}h_{1}(z_{j},z_{n})\\
\prod_{j\neq i,n}h_{1}(z_{j},z_{n})\prod_{j\neq i}h_{1}(z_{i},z_{j}) & =  \prod_{j\neq n}h_{1}(z_{j},z_{n})\prod_{j\neq i,n}h_{1}(z_{i},z_{j})
\end{align*}
which is always true. In the second case $w_{k}=z_{i}^{-1}$:
\begin{equation}
\frac{\prod_{j\neq i,n}h_{1}(z_{j},z_{n})}{\prod_{j\neq i,n}h_{1}(z_{i},z_{j})}\prod_{j\neq k}h_{q}(w_{j},z_{i})\prod_{j\neq i}h_{1}(z_{i}^{-1},z_{j})=\prod_{j\neq k}h_{q}(w_{j},z_{i}^{-1})\prod_{j\neq n}h_{1}(z_{i},z_{n})
\end{equation}
multiplying both sides by $z_{i}^{2(n-1)}$ and knowing that $z_{i}^{2}h_{1}(z_{i}^{-1},x)=h_{1}(z_{i},x)$ and $z_{i}^{2}h_{q}(x,z_{i}^{-1})=h_{q}(x,z_{i})$ we obtain the same equality.

Finally we calculate $g_n$:
\begin{equation}
 g_{n}=\frac{\prod_{j\neq n}h_{1}(z_{j},z_{n})\prod_{j}h_{q}(w_{j},z_{n})}{\prod_{j\neq n}h_{1}(z_{n},z_{j})\prod_{j}h_{q}(w_{j},z_{n})}=(-1)^{n-1}
\end{equation}
we replace $\sum_{i}g_{i}f(w_{j},z_{i})$ by $g_{n}f(w_{j},z_{n})$:
\begin{align}
B_{n} & =  \sum_{j}(-1)^{n+j}\frac{q^{n(n-1)/2}\mathbf{f}_{n-1,\hat{w}_{j}\hat{z}_{n}}f(w_{j},z_{n})}{\left(\prod_{i<k}h_{1}(z_{i},z_{k})(1-q^{2}w_{i}w_{k})\right)} \\
 & = q^{\frac{n(n-1)}{2}}\frac{\mathbf{f}_{n}}{\left(\prod_{i<k}h_{1}(z_{i},z_{k})(1-q^{2}w_{i}w_{k})\right)}\nonumber
\end{align}

\subsection{Integral version}
A special case (of direct interest to us) is when we integrate $B_n$ on a contour which surrounds only the poles $w_i=z_j^{\pm 1}$. Let us thus consider the following integral
\begin{equation}\label{intver}
 \oint\ldots\oint\prod_i  \frac{dw_i}{2\pi i} \psi(w,z) B_n(w,z) 
\end{equation}
where $\psi(w,z)$ is an analytic function of the $w$ in the integration region. Looking at the expression~\eqref{f}, we note that if in the calculation of $\mathbf{f}_n$ we pick a term with at least one $h_q(w_i,z_j)$ there will be fewer than $n$ poles and the integral will be zero. This way, we can erase all the terms with $h_q(w_i,z_j)$, and form the restricted $\mathbf{\bar{f}}_n$:
\[\mathbf{\bar{f}}_n=\frac{1}{(q^{-2}-1)^{n}\prod_{i}z_{i}(1-q^{2}w_{i}^{2})}\det\left|\frac{1}{h_{1}(w_{i},z_{j})}\right|\]
If we rewrite $h_{1}(w_{i},z_{j})=w_{i}z_{j}(w_{i}+w_{i}^{-1}-z_{j}-z_{j}^{-1})$ we easily identify $\mathbf{\bar f}$ with a Cauchy determinant, which can be evaluated:\begin{align*}
\mathbf{\bar{f}}_n & =  \frac{1}{(q^{-2}-1)^{n}\prod_{i}z_{i}^{2}w_{i}(1-q^{2}w_{i}^{2})}\frac{\prod_{i<j}(w_{i}+w_{i}^{-1}-w_{j}-w_{j}^{-1})(z_{j}+z_{j}^{-1}-z_{i}-z_{i}^{-1})}{\prod_{i,j}(w_{i}+w_{i}^{-1}-z_{j}-z_{j}^{-1})}\\
 & =  \frac{1}{(q^{-2}-1)^{n}\prod_{i}z_{i}(1-q^{2}w_{i}^{2})}\frac{\prod_{i<j}h_{1}(w_{i},w_{j})h_{1}(z_{j},z_{i})}{\prod_{i,j}h_{1}(w_{i},z_{j})}\end{align*}
Thus, in Eq.~\eqref{intver} one can rewrite $B_n$, given in general by Eq.~\eqref{Bf}, as
the same expression in which $\mathbf{f}_n$ is replaced with $\mathbf{\bar f}_n$.

Let us now assume that $\psi$ is of the form $\psi(w,z)=\prod_{i<j}(w_j-w_i)\,\phi(w,z)$
where $\phi$ is symmetric in the $w_i$. Then
\begin{multline}\label{inhoas}
\oint\ldots\oint\prod_i  \frac{dw_i}{2\pi i} \psi(w,z) B_n(w,z)\\
=n!
\oint\ldots\oint\prod_{i}\frac{dw_{i}}{2\pi i}\phi(w,z)\frac{\prod_{i<j}(w_{j}-w_{i})(qw_{i}-q^{-1}w_{j})}{\prod_{i\le j}h_{1}(w_{j},z_{i})\prod_{i\ge j}h_{q}(w_{j},z_{i})} \\ 
= \frac{1}{(q^{-2}-1)^{n}}\oint\ldots\oint\prod_{i}\frac{dw_{i}}{z_{i}2\pi i}
\phi(w,z)
\frac{q^{\frac{n(n-1)}{2}}\prod_{i<j}(w_{j}-w_{i})h_{1}(w_{j},w_{i})}{\prod_{i\leq j}(1-q^{2}w_{i}w_{j})\prod_{i,j}h_{1}(w_{i},z_{j})}
\end{multline}

\subsection{Homogeneous Limit}

The case of interest to us is when we set all the $z_i=1$. 
We can then use the same transformation as before:
\begin{equation*}
 u_i  =  \frac{w_i - 1}{q w_i -q^{-1}}
\end{equation*}
to deduce the desired equation from Eq.~\eqref{inhoas}. 

Call $\varphi(u)=\prod_i (1-q u_i)^2 \phi(w_i=\frac{1-q^{-1}u_i}{1-q\, u_i},z_i=1)$.
The second line becomes
\[
n!
\oint\ldots\oint\prod_{i}\frac{du_{i}}{2\pi i u_{i}^{2i}}\frac{\varphi(u)}{(q-q^{-1})^{n(n+2)}}\prod_{i<j}(u_{j}-u_{i})(1+\tau u_{j}+u_{i}u_{j})\]
while the expression on the third line becomes
\[
\oint\ldots\oint\prod_{i}\frac{du_{i}}{2\pi i u_{i}^{2n}}\frac{\varphi(u)}{(q-q^{-1})^{n(n+2)}}\frac{\prod_{i<j}(u_{j}-u_{i})(u_{i}-u_{j})(u_{i}+u_{i}+\tau u_{i}u_{j})}{\prod_{i\le j}(1-u_{i}u_{j})}\]
In both cases, the integrals surround zero.

In the latter, one can reinterpret some factors as a Vandermonde determinant:
\begin{equation*}
 AS\left\{ \prod_{i}\frac{(1+\tau u_{i})^{i-1}}{u_{i}^{2i}}\right\} =\prod_{i}\frac{1}{u_{i}^{2n}}\prod_{i<j}(u_{i}-u_{j})(u_{i}+u_{j}+\tau u_{i}u_{j})
\end{equation*}
and replace to obtain our final result:
\begin{multline}
\oint\ldots\oint\prod_{i}\frac{du_{i}}{2\pi i}\frac{\varphi(u)}{u_{i}^{2i}}\prod_{i<j}(u_{j}-u_{i})(1+\tau u_{j}+u_{i}u_{j})\\=\oint\ldots\oint\prod_{i}\frac{du_{i}}{2\pi i}\varphi(u)\frac{(1+\tau u_{i})^{i-1}}{u_{i}^{2i}}\frac{\prod_{i<j}(u_{j}-u_{i})}{\prod_{i\le j}(1-u_{i}u_{j})}
\end{multline}
where we recall that $\varphi(u)$ is some analytic function in a neighborhood of zero (that is, without poles in this domain) and symmetric in the $u_i$.

\bibliography{TDF}
\bibliographystyle{amsplain}

\end{document}